\documentclass[a4paper]{article}
\usepackage{authblk}
\usepackage{threeparttable}
\usepackage{graphicx}
\usepackage{lscape}
\usepackage{natbib}
\usepackage[english]{babel}
\usepackage[utf8]{inputenc}
\usepackage{lipsum}
\usepackage[T1]{fontenc}
\usepackage{amsmath}
\usepackage{graphicx}
\usepackage{subfigure}
\usepackage{amsmath}
\usepackage{amsfonts}
\usepackage{rotating}
\usepackage{rotfloat}
\usepackage{textcomp}
\usepackage{siunitx}
\usepackage{fontenc}
\usepackage{multirow}
\usepackage[colorinlistoftodos]{todonotes}
\usepackage{verbatim}
\usepackage{dcolumn}
\usepackage{rotating}
\usepackage{array}
\usepackage{booktabs,dcolumn}
\usepackage{setspace}
\usepackage{todonotes}

\newcolumntype{d}{D{.}{.}{2.5}}           

\newcommand{\BB} {{\boldsymbol \beta}}
\newcommand {\xx} {\mathbf x}
\newcommand {\XX} {\mathbf X}

\newcommand {\zz} {{\mathbf z}}
\newcommand {\yy} {{\mathbf y}}
\newcommand {\VV} {{\mathbf V}}
\newcommand {\ZZ} {{\mathbf Z}}
\newcommand {\II} {{\mathbf I}}
\newcommand {\WW} {{\mathbf W}}
\newcommand {\UU} {{\mathbf U}}
\newcommand {\DD} {{\mathbf D}}
\newcommand {\rr} {{\mathbf r}}
\newcommand {\gam} {{\boldsymbol \gamma}}
\newcommand {\Gam} {{\boldsymbol \Gamma}}
\newcommand {\su} {\sigma^2_\gamma}
\newcommand {\se} {\sigma^2_\varepsilon}
\newcommand {\suq} {\sigma^2_{\gamma_q}}
\newcommand {\seq} {\sigma^2_{\varepsilon_q}}
\newcommand {\suqt} {\sigma^{2(t)}_{\gamma_q}}
\newcommand {\seqt} {\sigma^{2(t)}_{\varepsilon_q}}
\newcommand {\SSigma} {{\boldsymbol \Sigma}}

\title{The use of sampling weights in the M-quantile random-effects regression: an application to PISA mathematics scores }

\author[$\textasteriskcentered$] {Francesco Schirripa Spagnolo}
\author[$\dag$] {Nicola Salvati}
\author[$\ddag$] {Antonella D'Agostino}
\author[$\S$] {Ides Nicaise}

\affil[$\textasteriskcentered$] {{\small Dipartimento di Studi Economici e Giuridici, Universit\`{a} di Napoli ``Parthenope''}}
\affil[$\dag$] {{\small Dipartimento di Economia e Management, Universit\`{a} di Pisa}}
\affil[$\ddag$] {{\small Dipartimento di Studi Aziendali e Quantitativi, Universit\`{a} di Napoli ``Parthenope''}}
\affil[$\S$] {{\small KU Leuven}}

\date{}
\begin{document}
	\maketitle

\begin{abstract}
M-quantile random-effects regression represents an interesting approach for modelling multilevel data when the interest of researchers is focused on the conditional quantiles. When data are based on complex survey designs, sampling weights have to be incorporate in the analysis. A pseudo-likelihood approach for accommodating sampling weights in the M-quantile random-effects regression is presented. The proposed methodology is applied to the Italian sample of the ``Program for International Student Assessment 2015'' survey in order to study the gender gap in mathematics at various quantiles of the conditional distribution. Findings offer a possible explanation of the low share of females in ``Science, Technology, Engineering and Mathematics'' sectors. 

\vspace{0.2cm}

\end{abstract}

\section{Background} \label{intro}

Gender differences in educational outcomes have raised concern over the last few decades because several studies stressed that the gender gap in education contributes to the gender segregation in the labour market and to the gender differences in wages \citep{brown1997, duquet2010, wb2011}. 

Especially, the severe female disadvantage in mathematics is of particular importance for its snowball effect. Although in many western countries the share of women that participate in higher education is now large, they yet tend to choose the STEM (Science, Technology, Engineering, and Mathematics) sectors much less frequently than men \citep{bradley2000, van2006, she2015}.
Consequently, the low share of women in STEM education is one of the reasons of the poor female presence in research and innovation decision-making positions \citep{van2006, she2015}. 
Since both research and innovation are important elements for a smarter, more sustainable and inclusive growth, the \cite{eu2012} also stressed that the female presence among scientific decision makers is fundamental for increasing the quality of research and improving the acceptance of innovation in the labour market. 

\begin{figure} [H]
	\centering
	{\includegraphics[width=0.6\textwidth]{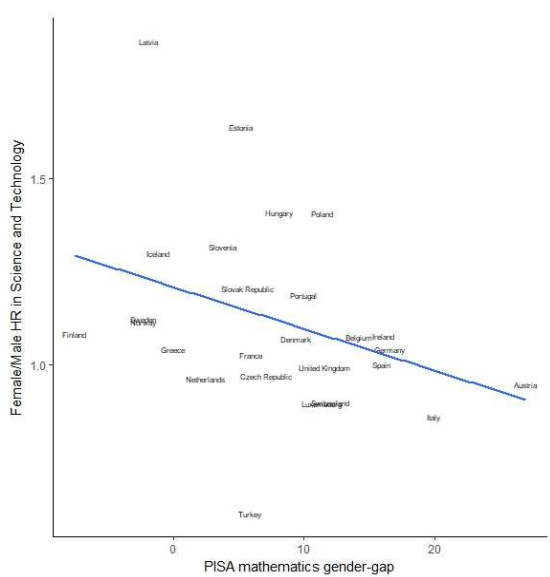}}
	\caption{Plot of gender gap in mathematics and ratio of female and male human resources in science and technology aged 25-64 in 26 OECD EU Countries (the data about the Human Resources in Science \& Technology refer to 2016. Source: Eurostat)\label{fig:gap}}
\end{figure}

Fig. \ref{fig:gap} highlights the negative relationship between the female-male ratio of people being employed in science and technology and the gap in mathematics in secondary school measured using data of the ``Program for International Student Assessment (PISA) 2015'' survey. As shown in Fig. \ref{fig:gap}, Italy is one of the countries with the highest gender gap in mathematics and the lowest share of female being employed in Science and Technology. 

Despite the gender differences in mathematics are much higher than the OECD average (a difference of 20 points against an OECD-average difference of 9 points) and the share of girls that does not achieve the minimum level of competencies in mathematics is one of highest among the OECD countries, in Italy the scores in mathematics are similar to the OECD average  (Italy has an average test score of 490, which is not significantly different from the OECD mean) \citep{invalsi2015}.

For this reason, the gender differences in mathematics have been studied in several papers: among others, \cite{bratti2007} used multilevel approach for analysing the 2003 Italian wave of the PISA and they found a large gender gap in mathematical performance;  \cite{agasisti14} found a female penalty on mathematics achievement stronger in Southern  than Northern regions of students attending grade 6 in the year 2011/12; \cite{masci2016} fitted a bivariate multilevel linear model for students attending the first year of junior secondary school in the year 2012/2013 and their results highlighted that on average males are better in mathematics, while females are better in reading; 
\cite{grilli2016} analysed 
the combined dataset of International Association for the Evaluation of Educational Achievement (IEA) assessments of student achievement in mathematics and science (TIMSS) and in reading (PIRLS) and they found that females have a lower performance in mathematics and science, but not in reading.

However, all these studies are also based on standard multilevel modelling approaches and accordingly   findings are valid mainly around the centre of the distribution  \citep{goldstein2011multilevel}. As a consequence, this can lead to an incomplete overview of the gender gap in mathematics, as large gender gap in mathematics among high-performing students can be more serious than among low-performers, since the best students tend to enroll more into higher education \citep{alma2017}.

The report of \cite{invalsi2015} showed, based on simple raw-data, how the gender gap in mathematics really tends to be wider in the upper quantiles of the distribution. In this perspective, what is currently lacking in the literature is the analysis of the gender gap along the overall conditional distribution of mathematics scores, i.e. taking into account other relevant variables (e.g. students and school-level characteristics) that can mitigate the gender effect in different quantiles. 

In other words the aim of this paper is to estimate and understand the gender gap in math, using a distributional approach, which makes comparisons over the entire range of mathematics scores rather than focusing on summary measures such as the mean. Such analysis is obviously relevant as it might give insights to policy makers for implementing policies to reduce the gap and to try and raise the number of girls enrolled in scientific graduate programmes.

From a methodological point of view, quantile regression is the most popular approach for estimating the conditional quantiles of a response variable  \citep{koenker1978, koenker2005}, but the classical implementation of this estimation method does not allow to include into the analysis specific random effects to take into account the hierarchical data structure. Actually, some attempts to implement quantile regression with a hierarchical data structure have been addressed in literature. For instance, \cite{koenker2004} proposed a penalized fixed-effects estimation method for longitudinal data, while \cite{geraci2007, geraci2014} proposed a two-levels quantile regression with mixed effects named the Linear Quantile Mixed Model (LQMM).  It allows to model data with complex dependence structures by including multiple random effects in the linear conditional quantile functions. For a review of linear quantile models for multilevel data see \citep{marino2015}. In the context of educational research, there are not many examples of studies that applied multilevel quantile regression models. \cite{costanzo2015} used the LQMM to evaluate the effect of a specific training programme (M@tabel) on the Italian sixth grade students' performance in mathematics at secondary schools: the author highlighted the advantages of using this model compared to the traditional linear random effects one. \cite{faria2016} applied LQMM to analyse the determinants of students’ success in Portugal using 2009 PISA survey.

A more flexible alternative to quantile regression is also the the M-quantile (MQ) regression \citep{breckChamb88, ChambTzavidis2006Mq}. This approach integrates the concepts of quantile regression and expectile regression \citep{Newey1987} within a framework defined by a `quantile-like' generalization of regression based on influence functions (M-regression).
The classical MQ regression even does not allow the analysis of multilevel data. Recently important developments have been made in order to overcome this limitation. \cite{tzavidis2016} extended M-quantile regression so as to include random effects in two-level hierarchical data structure using maximum likelihood (MQRE-2L); \cite{borgoni2016} extended the MQRE-2L to a three-level random effects model (MQRE-3L); \citet{alfo:2017} define a finite mixture of quantile and M-quantile regression models (FMMQ) for heterogeneous and /or for dependent/clustered data. Although LQMM, MQRE and FMMQ allow for modelling data with a complex dependence structure  including random effects, they do not accommodate for sampling weights in the estimation procedure. Thus, they are not suitable for modelling survey data based on complex survey designs (e.g. when statistical units are selected with unequal probabilities from the population) usually implemented by national and international organizations. 
Namely, any estimation procedure that does not take into account sampling weights provides biased results \citep{cochran1977, sarndal1992, pfeffermann1993}.

While weighting for unequal probabilities of selection is a relatively well-established procedure in single level models as it can be viewed as an application of the `Pseudo Maximum Likelihood' (PML) approach \citep{chambers2003,binder1983}, it needs additional remarks in multilevel models because the covariance structure of the population has to be modelled. For this reason, the multilevel models require the knowledge of the inclusion probabilities at each hierarchical level. \cite{pfeffermann1998weights} used the PML approach and the iterative generalized least square (IGLS) algorithm \citep{goldstein1986} in a two-level model with a continuous response variable, obtaining the so-called probability-weighted IGLS (PWIGLS) estimators. 

Alternative approaches have been proposed by \cite{grillipratesi, rhs2006, aspa2006}. These approaches include sampling weights in the estimation procedure. In particular, the census likelihood is estimated by weighting the sample likelihood including sampling weights in the log-likelihood function. Indeed, \cite{rhs2006} proposed a pseudo-likelihood approach via adaptive quadrature for generalized linear models with dichotomous response with any number of levels; their method also allows for stratification and PSUs that are not represented by a random effect in the model. \cite{grillipratesi} used the PML approach to develop weighted estimators in the context of ordinal and binary models.
Finally, \cite{aspa2006} proposed an approximately unbiased multi-level pseudo maximum likelihood (MPML) for general multi-level modelling and asserted that it can be used with any parametric family of distributions and any linear or non-linear multilevel models.   

The proposed multilevel approaches based on the use of the sampling weights suffer from the same limitations as the the standard multilevel models: they allow for modelling the centre of the distribution. This limitation can lead to an incomplete overview of the gender gap in mathematics. Therefore, the purpose of this paper is to extend the current modelling framework of M-quantile random effects so that to include sampling weights and to examine the shape of gender gap along the overall conditional distribution of mathematics scores using PISA-OECD 2015 data. 

The structure of the paper is as follows: after this introduction, Section \ref{data} provides details of PISA data, including variables used in the empirical analysis. Section \ref{methodology} provides  details of the methodology proposed. In Section \ref{results} we present the results of applying the new methodology to PISA data for studying gender gap in mathematics, together with details of the method that is used for model evaluation. In Section \ref{simulation} we empirically evaluate the properties of the proposed model by using a Monte Carlo simulation study. Finally, Section \ref{conclusion} provides a concluding summary and a discussion of potential areas for future research.

\section{Data and sample characteristics} \label{data}
The empirical analyses are based on Programme for International Students Assessment (PISA) survey that is a triennial survey started in 2000 and conducted by the OECD for measuring the extent to which the 15-year-old students have acquired key knowledge and skills in three school subjects (mathematics, science and reading) that ``\textit{are essential for full participation in modern societies}'' \citep{oecd2016} in many countries and economies. In particular, our data refer to the latest 2015 PISA wave.

PISA uses a two-stage stratified sample design. Schools having 15-year-old students are the first-stage sampling units. A systematic Probability Proportional to Size (PPS) sampling is used to select schools. In order to improve the precision of sample-based estimates, schools previously are divided into explicit strata based on school characteristics. Students within sampled schools are the second-stage sampling units. Although the two-stage sampling design used in PISA should guarantee that the students have the same probability of selection, three factors contribute in explaining the variability of weights: i) over- or under-sampling of some strata of the population; ii) lack of accuracy or no update in defining the size of school; iii) adjustment of weights for school and student non response. Further details on sample design and weights defined in PISA can be found in \citet[cap. 3]{pisaspss}; \citet[cap. 4 and 8]{tech2015}; \cite{rhs2006}.

Although the target population of PISA survey is represented by 15-years-old students who have completed at least 6 years of formal schooling, our analysis is based only on students enrolled in public upper secondary schools, corresponding to the 3th level of the International Standard Classification of Education (ISCED). This choice is determined by the fact that upper secondary education represents a sort of minimum credential for entering successfully into the labour market and it is required  to pursue their university education. For this reason, the upper secondary level of education is considered a crucial indicator of the output of the educational system of a Country  \citep{oecd2003}. Our estimation sample contains 7163  students with non-missing response variable and covariate information. The 7163 students are clustered in 283 schools.

The response variable that we use in this paper is the students' scores in mathematics. In particular, we use the first plausible value generated by the standardized procedure implemented by the PISA team (for more detail about the methodology see \citealp{tech2015}).

The covariates used in the models are as follows. Gender (male (baseline) or female) whose effect is the main object of this study. Immigration status was coded as a dummy variable \citep{oecd2016}. Namely, non immigrant students (students whose mother or father (or both) were born in the country where the PISA test occurs, regardless of whether the student himself or herself was born in that country -- baseline) \textit{vs} immigrant students (students whose mother and father were both born in a country other than that where the student performs the PISA test). The native-immigrant gap is well recognized in the literature (see for instance \citealp{rangvid2007, ammermueller2007, azzolini2012, hajisoteriou16}) because usually immigrants encounter many difficulties in achieving adequate mathematical skills, consequently they cannot participate successfully and actively in the host society. The index of schoolwork-related anxiety was measured by the responses of the students to the statements regarding their worry about tests and study in general. It was standardised to have a mean of 0 and a standard deviation of 1 across OECD countries. Positive values on the index indicate that students reported higher levels of schoolwork-related anxiety than the average student across OECD countries; negative values indicate that students reported lower levels of anxiety than the average student \citep{oecd2016volIII}. Grade repetition is a good indicator of the students' school career and it can affect negatively the academic performance and students' delayed entry into the labour market  (\citealp{oecd2016volII, oecd2012IV}; see also \cite{miyako2014} for a review of studies about this topic). It is recoded into two categories: the student has never repeated a grade in any level (baseline) and the student has repeated a grade in at least one level. Lack of punctuality was measured by students' answers on whether they had arrived late for school in the two weeks before the test (no (baseline) or yes). This indicator measures the students' truancy and, at the same time, it may represent the lack of interest for learning, moreover it has negative consequences on other students because it can contribute to a disruptive learning environment. 
Student socio-economic status, measured by the PISA index of economic, social and cultural status (ESCS). It was derived from the combination of three variables related to family background: highest parental education, highest parental occupation, and home possessions (as a proxy of family wealth). 
Higher values of ESCS indicate better socio-economic status \citep{oecd2016volIII}. 
This variable is strongly associated with the students' achievement because students with a high ESCS have better access to educational resources provided by their family (material resources and the educational level of their parents). 

In addition to individual-level characteristics, the key explanatory variables at school-level are as follows. The macro region  that indicates where the school is located (Northern regions (baseline) or Southern regions). Regional disparities in school achievement are indeed, generally very serious in Italy:  students attending the schools located in Northern Italy tend to have higher achievement scores than their counterparts in Southern Italy \citep{invalsi2012, invalsi2015}. The type of school is recoded into three categories: Lyceums (baseline), Technical and Vocational Schools and Other Schools. This classification allows to discriminate between three groups of schools that offer different potential outcomes. Indeed, Lyceums provide mostly a theoretical training and prepare students for tertiary education (ISCED 5 and 6); Technical and Vocational Schools offer both a general education and a technical specialization in a specific field of study; other schools are institutions designed to offer technical activities and allow students to obtain a professional qualification which is immediately recognised in the labour market.

Moreover, the school-mean of each individual-level variable has been calculated for estimating the so called ``compositional effect'' \citep{sani2011}. These variables are indeed proxies of the social context of the school that might considerably influence the academic achievements because of peer interactions. Moreover, on the basis of school composition teachers can change their instructional practices to take into account the characteristics of school and consequently this phenomenon can affect the individual students' performance \citep{rangvid2007quantile,schneeweis2007}. 
Therefore, how the social context of the school affects the achievement of pupils is a relevant element for policy makers in designing measures  to control the issue of social composition of the schools. Table \ref{sum1} describes the sample by reporting sample proportions for binary variables, means, medians and quantiles for continuous variables.

\begin{table}\caption{\label{sum1}Descriptive statistics of outcome variable and covariates} 
	\centering
	\footnotesize
\begin{tabular}{lrrrrrr}
					 \toprule
					 \noalign{\smallskip}
					 \multicolumn {1}{l}{Variable} & \multicolumn {1}{c}{Min} & \multicolumn {1}{c}{$1^{st}$Q}& \multicolumn {1}{c}{Median}& \multicolumn {1}{c}{Mean}& \multicolumn {1}{c}{$3^{rd}$Q}& \multicolumn {1}{c}{Max.}\\	
					 \midrule \noalign{\smallskip}	
					 MATHEMATICS    &140.80&432.14&498.17&496.34&560.72&822.64\\
					 ESCS	&-2.99&-0.71&-0.04&-0.05&0.67&3.56\\
					 Anxiety index	&-2.51&-0.08&0.52&0.48&1.05&2.55\\
					 Female  &&&&0.51&&\\
					 Immigrants &&&&0.07&& \\
					 Grade repetition  &&&&0.13&& \\
					 Lack of punctuality  &&&&0.35&& \\		
					 South Italy    &&&& 0.44&&\\
					 \multicolumn {1}{l}{Type of School} &&&& &&\\
					 \hspace{0.3cm}Lyceums &&&& 0.42&&\\
					 \hspace{0.3cm}Technical and Vocational	&&&& 0.54&&\\
					 \hspace{0.3cm}Others &&&& 0.04&&\\
					 Mean gender &&&& 0.50 &&\\
					 Mean Immigrant status &&&& 0.09&&\\
					 Mean Grade repetition &&&& 0.16&&\\
					 Mean Lack of punctuality &&&& 0.36&&\\
					 Mean ESCS &-1.50&-0.59&-0.27& -0.20&0.15&1.07\\
					 Number of units & 7163 &&&&&\\
					 Number of clusters & 283 &&&&&\\
					 \bottomrule 
\end{tabular}
\end{table}

There is an approximately even split of male and female students, only the 7\% of the students' population has a foreign background  in line with the foreign citizens who live in Italy (cfr. \citealp{istat2017}), a little bit more than  a third of students reported to arrive late at school and only the 13\% of them have repeated a grade. The average score in mathematics is almost 450 and about 75\% of all students exceed the level 2 of proficiency (higher than 420.07 and less than or equal to 482.38), in line with the OECD mean (cfr. \citealp{invalsi2015}). The inter-quantile range is 129 a little more than OECD mean difference (125 points). More than 50\% of all schools are located in Northern/Central regions and more than half (54\%) are Technical and Vocational Schools.

Finally, to support the choice of using a robust approach a simple weighted two-level random intercept model for mathematics score has been estimated using the student-level and the school-level weights provided in PISA database and the covariates described above.

\begin{figure} [H]
	\subfigure[\label{fig:resmath}]%
	{\includegraphics[width=1\textwidth]{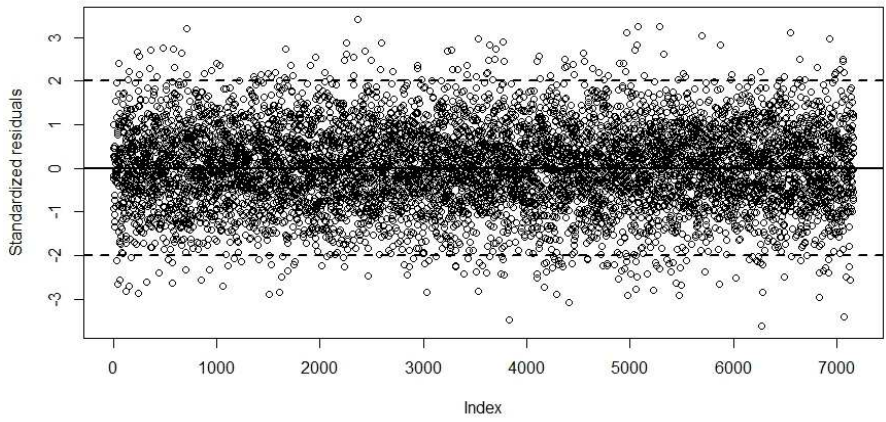}}
	\subfigure[\label{fig:cookmath1}]%
	{\includegraphics[width=1\textwidth]{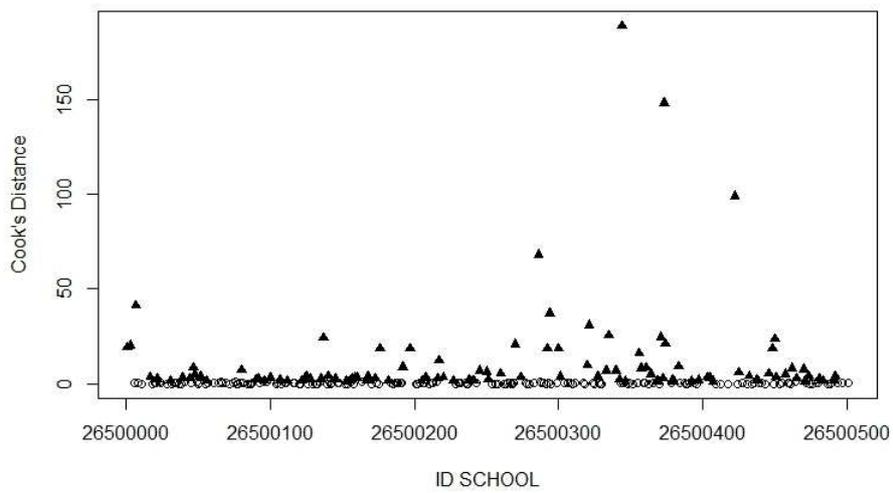}}
	\caption{Plot of standardised residuals (a) and Cook's Distance (b) derived by fitting a linear mixed-effects model for mathematics scores\label{fig:math}}
\end{figure}

Fig. \ref{fig:math} shows the residual analysis. Clearly, several outliers and second-level units that can be classified as influential points (filled triangle points exceed the cut-off equal to 1) stand out.

\section{Methodology} \label{methodology}
In this section, after a briefly reviewing of the PML in multilevel models (Section \ref{sec:PML}), we present the weighted two-level MQRE model (Section \ref{W_MQRE}).
\subsection{Pseudo-maximum likelihood approach in multilevel model} \label{sec:PML}
Let us have a two-level population with $M$ level 2 units (or cluster) and $N_j$ level 1 units within the $j${th} cluster. 
Suppose that an outcome continuous variable $y$ is modelled 
using a two-level model:
\begin{equation} \label{LMM}
y_{ij}=\xx^T_{ij}\BB + \zz^T_{ij}\gam_j + \varepsilon_{ij}, \qquad i=1,...,N_j, \qquad j=1,...,M,
\end{equation}
where $\xx^T_{ij}$ is a covariate-vector of dimension $p$; $\BB$ is a vector of regression coefficients; $\zz_{ij}$ is a vector of group indicators used to define the random part of the model; $\gam_j$ are random effects varying over clusters and we assume that they follow a multivariate normal distribution with zero mean and covariance matrix $\mathbf{\Omega}$ and $\varepsilon_{ij}$ are individual random effects and we assume $\varepsilon_{ij} \sim N(0, \sigma_{\varepsilon}^2)$. 
In the paper we consider a two-level random-intercept model for unit $i$ in cluster $j$. It can be written as
\begin{equation} \label{LMM_ri}
y_{ij}=\xx^T_{ij}\BB + \gamma_j + \varepsilon_{ij},
\end{equation}
where  $\gamma_j \sim N(0,\sigma^2_\gamma)$.
When population model in Eq. \ref{LMM_ri} is evaluated at the sample, estimation of the parameters of the model can be obtained by employing maximum likelihood method. Under the normal assumption and $\varepsilon_{ij}$ independent from $\gamma_j$ the log-likelihood function is:
\begin{equation} \label{logLMM}
log L(\BB,\su, \se)= - \frac{1}{2}\mathrm{log}|\VV|-\frac{1}{2}(\yy-\XX\BB)^{T} \VV^{-1}(\yy-\XX\BB), 
\end{equation}
where $\yy$ i the $n \times 1$ response vector and $n$ is the sample size; $\VV= \SSigma_\varepsilon + \ZZ \SSigma_{\gamma} \ZZ^T$, $\SSigma_\varepsilon= \se\II_{n}$, $\SSigma_\gamma= \su\II_{M}$ and $\ZZ$ is an $n \times M$ matrix of known positive constant. Here $\II_{k}$ represents the identity matrix of dimension $k$. Estimates of fixed and variance parameters can be obtained by differentiating the log-likelihood with respect to these parameters and then solving the estimating equations defined by setting these derivatives equal to zero.  
However, if the sample data are generated by a complex sample design the different inclusion probabilities have to take into account and the weights have to be inserted somewhere in the estimation process.

This idea leads to use the weights within regression analysis to estimate the model on the entire population (a census) rather than the sample. The PML approach in multilevel framework requires to write down the census-likelihood function and then incorporates weights at each level of the analysis in the log-likelihood equations. 

The census-likelihood function can not be expressed as a simple sum of the elementary unit contributions, but it is a sum across all levels of the data hierarchy. In a two level model it is a function of sums across level 2 and level 1 units:
\begin{equation} 
\log L(\BB, \su, \se)= \sum_{j=1}^{M} \log\int \left[\exp \left\{\sum_{i=1}^{N_{j}} \log L_{ij}(\BB, \su, \se|\gam_j)\right\} \right] f(\gam_j)\mathrm{d}\gam_j,
\end{equation}
where $f(\gam_j)$ is the normal probability density of the level 2 random effects and \linebreak $\log L_{ij}(\BB, \su, \se|\gam_j) $ is the log-likelihood contribution of the level 1 units conditioned on the level 2 random effects. 

Now, suppose that the whole population is not observed, but the sample data have been selected with the following two-stage sampling design: at the first stage $m$ clusters are selected with inclusion probabilities $\pi_j$, $j=1,...M$, and at the sub-sequent stage elementary units are sampled with conditional probabilities $\pi_{i|j}$, $i=1,...N_{j}$. The log-likelihood function for the sample units can be expressed as  
\begin{equation} \label{wL}
\log L_w(\BB, \su, \se)= \sum_{j \in s} w_j \log\int \left[\exp \left\{\sum_{i \in s_j} w_{i|j} \log L_{ij}(\BB, \su, \se|\gam_j)\right\} \right] f(\gam_j)\mathrm{d}\gam_j,
\end{equation} 	
where $s$ and $s_j$ indicate the sampled clusters and the units sampled in the cluster $j$, respectively. 
As we noted, the sum for each of the levels requires the respective
conditional probabilities of selection: the log-likelihood contributions of the level-1 units are weighted by $w_{i|j}= 1/\pi_{i|j} $, i.e. the inverse of the selection probability of unit $i$ in cluster $j$ given that cluster
$j$ has been sampled, and the log-likelihood contributions of the level 2 units are weighted by $w_j=1/\pi_j$, i.e. the inverse of the selection probabilities of cluster $j$.

In case of two-level random-intercept model the Eq. \ref{wL} can be expressed as:
\begin{equation} \label{logw}
\log L_w(\BB, \su, \se) = -\frac{1}{2} \sum_{j \in s} w_j \left[ \log(2\pi |\WW^{-1} \VV |) + (\yy-\XX\BB)^{T} \VV^{-1}(\yy-\XX\BB) \right],
\end{equation}
where $\WW$ is a diagonal matrix of the first-level sampling weights, $\WW =\mathrm{diag}(w_{1|j},$\\$..., w_{n_j|j})$, and $\VV=\WW^{-1}\Sigma_\varepsilon + \ZZ \Sigma_{\gamma} \ZZ^T$.

Differentiating Eq. \ref{logw} with respect to the fixed effects and to the variance components the weighted estimation equations of the LMM are obtained. 

\subsection{Weighted \textit{M}-quantile Random Effects Model} \label{W_MQRE}
As we explained in Section \ref{intro}, the classical multilevel regression model in Eq. \ref{LMM} provides an incomplete picture of the distribution of the response variable given the auxiliary information, because it just summarises the behaviour of the mean of the outcome variable at each point in a set of covariates. Moreover, the presence of outliers in the data invalidates the general assumptions of the model; in this case the estimators of the parameters of the model under Eq. \ref{logLMM} could be biased and inefficient \citep{richardsonwelsh1995}. \cite{huggins1993} and \cite{richardsonwelsh1995} proposed an approach based on M-estimation for a robust estimation of the multilevel models. \cite{tzavidis2016} used this idea of robust estimation in the context of multilevel model to extend M-quantile regression in order to include random effects that account for a two-level hierarchical structure in the data and used maximum likelihood technique to estimate the parameters of the model.

In the linear case the M-quantile regression assumes the following form:
\begin{equation}  \label{mq}
MQ_y(q|\xx;\psi)=\xx_i^T\BB_{\psi q},
\end{equation}
where $\BB_{\psi q}$ is defined as the minimiser of

\begin{equation}  \label{mq1}
\min_{\BB}E\left[|q-I(u<0)|\rho(u)\right],
\end{equation}
with $u=(y-\mathbf{x}^T\BB)/\sigma_q$, $\sigma_q$ a scale parameter. In this paper we use the popular Huber loss function \citep{breckChamb88}:
\begin{equation}\label{rho}
\rho(u)= \left\{
\begin{array}{ll}
2c|u|-c^2 & \quad |u|>c, \\
u^2   & \quad|u|\leq c, 
\end{array}
\right.
\end{equation}
where $c$ is a tuning constant bounded away from zero (a common choice is c = 1.345).
Assuming that $\rho$ is (a.e.) continuously differentiable and convex, an estimator of $\BB_{\psi q}$, $\hat \BB_{\psi q}$ can be obtained as the solution of the following system of equations
\begin{equation}
\sum_{i=1}^{n} \psi_q\left(\frac{y_i - \xx_i^T\BB_{\psi q}}{\hat{\sigma}_q}\right),
\end{equation}
where $\hat{\sigma}_q$ is a consistent estimator of ${\sigma}_q$, $\psi_q(u)$ denotes an asymmetric influence function, which is the derivative of an asymmetric loss function $\rho_q(u)=|q-I(u<0)|\rho(u)$.

Defining the loss function as in Eq. \ref{rho}, MQ regression model estimation is obtained weighting positive residuals by $q$ and negative residuals by $1-q$. The Iterative Weighted Least Square (IWLS) algorithm used to fit an M-quantile regression model guarantees convergence to a unique solution when a continuous monotone influence function (such as Huber function with $c>0$) is used \citep{kokic1997}. An other interesting benefit is the possibility to trade efficiency for robustness by setting the value of $c$. For example, if $c$ is close to zero, the robustness of the model increases and its efficiency decreases, then MQ moves towards quantile regression. On the other hand, setting $c$ large, the robustness of the model decreases and its efficiency increases, then MQ moves towards expectile regression.

For modelling multilevel structured data \cite{tzavidis2016} extended the linear specification of the model in Eq. \ref{mq} to allow for the inclusion of random effects (MQRE-2L) to account for a two-levels hierarchical structure in the data. The M-quantile random intercept model is defined as follows:
\begin{equation} \label{MQRE-2L}
MQ_{y_{ij}}(q|\xx_{ij}, \gamma_j; \psi)=\xx^T_{ij}\BB_{\psi q} + \gamma_j,
\end{equation}
where $\gamma_j$ is the random effect for cluster $j$. For obtaining the estimation equations for the regression coefficients and the variance parameters the authors extended the idea of asymmetric weighting of residuals by changing the estimating equation of the ML proposal II by \cite{richardsonwelsh1995} following \cite{sinha2009}. 

\cite{borgoni2016} extended the MQRE-2L to a three-level M-quantile random effects regression (MQRE-3L). The main difference with the model proposed by \cite{tzavidis2016} is an additional equation for estimating the further level and a more complex variance-covariance matrix of the outcome variable.

The aim of this paper is to extend the idea of \cite{tzavidis2016} and incorporate sampling weights in the model. So, following the pseudo-likelihood approach we propose an estimation procedure that adjusts for the effect of an informative sampling design on estimation in MQRE-2L.

In particular our approach follows the works of \cite{grillipratesi}, \cite{aspa2006} and \cite{rhs2006}. In particular, the sampling weights are inserted before the derivatives of the log-likelihood function are taken. Consequently, firstly, the weighted log-likelihood in Eq. \ref{logw} is differentiated obtaining the weighted estimation equation for a the parameter of the random intercept model and then, following the idea of \cite{tzavidis2016}, we propose a robustification of the weighted estimation equations for estimating the regression coefficients and the variance parameters obtaining in the previous step as follows:

\begin{equation} \label{beta}
\sum_{j \in s} w_j \left[ \XX_j^{T} \VV_{qj}^{-1} \UU_{qj}^{1/2} \psi_q (\rr_{qj})  \right] = \mathbf{0},
\end{equation}
\begin{equation*}
-\frac{1}{2} \sum_{j \in s} w_j \left[K_{2qj} \mathrm{tr}(\VV_{qj}^{-1} \ZZ_j\ZZ_j^T)-  \psi_{q} (\rr_{qj})^{T} \UU_{qj}^{1/2} \VV_{qj}^{-1} \ZZ_j\ZZ_j^T \VV_{qj}^{-1} \UU_{qj}^{1/2} \psi_ {q} (\rr_{qj}) \right]= \mathbf{0},
\end{equation*}
\begin{equation} \label{var}
-\frac{1}{2} \sum_{j \in s} w_j \left[K_{2qj} \mathrm{tr}(\VV_{qj}^{-1}  \WW_j^{-1})-  \psi_{q} (\rr_{qj})^{T} \UU_{qj}^{1/2} \VV_{qj}^{-1} \WW_{j}^{-1} \UU_{qj}^{1/2} \VV_{qj}^{-1} \psi_{q} (\rr_{qj}) \right] = \mathbf{0},
\end{equation}
where $\XX_j$ is the matrix of covariates in the cluster $j$; $\rr_{qj}= \UU_{qj}^{-1/2}(\yy_j-\XX_j\BB_{\psi q})$ is a vector of scaled residuals with components $r_{ijq}$; $\UU_{qj}$ is a diagonal matrix with diagonal elements equal to the diagonal elements of the variance-covariance matrix of the $j${th} cluster $\VV_{qj}= \WW_j^{-1}\SSigma_{\varepsilon_q} + \ZZ_j \suq \ZZ_j^T $;  $\SSigma_{\varepsilon_q}= \seq\II_{n_j}$, $\suq$ and $\seq$ are the M-quantile specific variance parameters; $\ZZ_j$ is an $n_j \times 1$ vector of known positive constants; $K_{2qj}= \mathrm{E}[\psi_{q} (e_j)\psi_{q} (e_j)^T]$ with $ e_j \sim N(\mathbf{0}, \II_{n_j})$. Eq. \ref{beta} and \ref{var} are the estimating equations of the Weighted M-quantile Random Effects Model (henceforth Weighted-MQRE). To obtain estimators of $\BB_{\psi q}$, $\suq$ and $\seq$ Eq. \ref{beta} and \ref{var} are solved iteratively. For Eq. \ref{beta} a Newton-Raphson algorithm is used and for Eq. \ref{var} the fixed-point iterative method is implemented for getting the estimates. Details about the algorithm can be found in the Appendix A. Inference for the model parameters is performed using a sandwich estimator following the proposal by \cite{tzavidis2016} (for more details see Appendix B).

As pointed out in \cite{tzavidis2016} and \cite{borgoni2016} and discussed in detail in \cite{jones1994}, we can consider M-quantiles equivalent to quantiles because both target the same part of the distribution. So, the estimates of  $\BB_{\psi q}$ can be interpreted how the effect of a unit change in $x$ on a given $q$-quantile of the distribution of $y$. 

\subsection{Design consistency of the Weighted-MQRE regression coefficients}

In this section we prove the design consistency of $\hat\BB_{\psi q}$ for the regression coefficients $\BB_{\psi q}$ of the Weighted-MQRE. The asymptotic design-based properties of  the regression coefficients in the M-quantile framework have been introduced by \cite{fabrizi2014w} and the consistency of estimators in the multilevel modelling has been demonstrated by \cite{pfeffermann1998weights}. From a theoretical point of view, the establishment of consistency properties in multilevel modelling requires that the number of groups, $m$, and the units in each group, $n_j$, increase; however as pointed out in \cite{pfeffermann1998weights} in practice $n_j$ are often small and the consistency of $\hat\BB_{\psi q}$ can be established when only $m$ increases and tends to $M$.
Following the work of \cite{fabrizi2014w} and \cite{wang2011} we define, for any $q \in (0,1)$, the population parameter $\mathbf{B}_{\psi q}$ as: 
\begin{equation}
\mathbf{B}_{\psi q} = inf\{\BB_{\psi q} : S_M(\BB_{\psi q}) \geq 0\},
\end{equation} 
where $S_M(\BB_{\psi q}) = M^{-1} \sum_{j=1}^{M} \XX_j \VV_j^{-1/2} \psi_q(\VV_j^{-1} (\mathbf{y}_j  - \xx_j^T \BB_{\psi q}))$. 

The population parameter $\mathbf{B}_{\psi q}$ can be estimated by 
\begin{equation}
\hat{\BB}_{\psi q}= inf\{\BB_{\psi q} : \hat{S}_M(\BB_{\psi q}) \geq 0\},
\end{equation}
with $\hat{S}_M(\BB_{\psi q}) = M^{-1} \sum_{j=1}^{m} w_j \XX_j \VV_j^{-1/2} \psi_q(\VV_j^{-1} (\mathbf{y}_j  - \xx_j^T \BB_{\psi q}))$. 

Assuming the technical conditions A1-A5 (described in Appendix C), for each $q \in (0,1)$ and for large closed interval $\Theta \in \mathbb{R}^p$, $\mathrm{sup}_{\BB_{\psi q} \in \Theta} |\hat{S}_M(\BB_{\psi q})- S_M(\BB_{\psi q})|= o_p(1)$, the estimator $\hat{\BB}_{\psi q}^{\mathrm{W-MQRE}}$ obtained solved the Eq. \ref{beta} is design $\sqrt{m}$-consistent for $\mathbf{B}_{\psi q}$, in the sense that $M^{-1}(\hat{\BB}_{\psi q}^{\mathrm{W-MQRE}}-\mathbf{B}_{\psi q})=O_p(m^{-1/2})$. 

Proof of this statement follows the results in \cite{wang2011} and \cite{fabrizi2014w}. 

\section{Application: Modelling the conditional distribution of mathematics scores} \label{results}
\subsection{Modelling approach}
In this section we apply the method proposed in the Section \ref{W_MQRE} to study the effect of gender along the entire distribution of mathematics scores using PISA-OECD 2015 data. In particular, taking into account the socio-economic issues and the literature review that were described in Section 1, we are mainly interested in measuring the gender gap in mathematics, controlling for individual and school-level characteristics.  In particular, adding the cluster mean of the individual variables we can estimate the so-called within-effect. A specific function in the statistical programming environment R \citep{Rcore} has been written and it is available from the authors upon request.

The following scaled version of first-level weights named “scaling method 2” in \cite{pfeffermann1998weights} has been used in order to fit mathematics scores: 

$$
w_{i|j}^{s2}= n_j w_{i|j} \Big( \sum_{i}^{n_j} w_{i|j}\Big)^{-1},
$$
where $n_j$ is the number of sample units in the $j$th cluster and  $w_{i|j}$ are the first-level weights.
\cite{pfeffermann1998weights} have proposed an alternative scaling method, named ``scaling method 1'': $w_{i|j}^{s1}= w_{i|j}  (\sum_{i}^{n_j} w_{i|j}) (\sum_{i}^{n_j} w_{i|j}^2)^{-1}$. It is worth noting that the results using this alternative scaling method are almost identical to those using ``scaling method 2''. However, simulation results in \cite{pfeffermann1998weights} suggest that ``scaling method 2'' works better than ``scaling method 1'' for informative weights. Moreover, the student-level weights provided in PISA (the variable called W\_FSTUWT) are student-level overall inclusion weights ($w_{ij}$ not $w_{i|j}$) adjusted for non-inclusion and non-participation of students. Consequently scaling the weights is a procedure to overcome the problem of not having the conditional sampling weights. Indeed, if student-level weights are rescaled, the model estimates are equivalent to those obtained when $w_{i|j}$ is available; the adjustment factors do not affect the rescaled version of the student-level weights \citep{rhs2006}.

\subsection{Estimation Results}

The estimates from the M-quantile random effect model with sampling weights are presented in Table \ref{MQRE-math}. 
The estimated coefficients are significantly different from 0 at each quantile with a different intensity for almost all the covariates; hence they indicate the usefulness of the estimation method used for the analysis.  

\begin{table}\caption{\label{MQRE-math}Results--Weighted-MQRE model for mathematics scores\dag}
	\centering
	\footnotesize
		\begin{tabular}{l ddddd}
			\toprule 
			&\multicolumn{5}{c}{Results for the following value of q:}\\
			\midrule
			Variable    & 0.1 
			&  0.25 & 0.5 & 0.75 & 0.9\\
			\midrule
			Intercept	&	533.20	^{***}	&	570.63	^{***}	&	608.38	^{***}	&	648.20		^{***}	&	688.40		^{***}			\\
			Gender	&	-14.38	^{***}	&	-14.98	^{***}	&	-17.08	^{***}	&	-17.90		^{***}	&	-16.95		^{***}			\\
			Immigrant status	&	-14.58	^{**}	&	-15.53	^{**}	&	-14.95	^{***}	&	-15.19		^{***}	&	-18.74		^{***}			\\
			Grade repetition	&	-45.21	^{***}	&	-45.21	^{***}	&	-42.24	^{***}	&	-38.31		^{***}	&	-35.79		^{***}			\\
			Anxiety index	&	-5.29	^{**}	&	-5.82	^{***}	&	-6.74	^{***}	&	-8.20		^{***}	&	-9.66		^{***}			\\
			Lack of punctuality	&	-12.61	^{***}	&	-9.89	^{***}	&	-9.17	^{***}	&	-13.82		^{***}	&	-16.64		^{***}			\\
			ESCS	&	2.51		&	2.39		&	2.39		&	3.72			&	6.07		^{**}			\\
			Mean gender	&	-48.93	^{***}	&	-51.64	^{***}	&	-53.48	^{***}	&	-58.47		^{***}	&	-66.52		^{***}			\\
			Mean immigrant status	&	-137.76	^{***}	&	-114.99	^{**}	&	-81.98	^{***}	&	-78.28		^{***}	&	-80.52		^{***}			\\
			Mean grade repetition	&	-51.42	^{*}	&	-54.71	^{*}	&	-62.36	^{**}	&	-63.36		^{**}	&	-67.51		^{***}			\\
			Mean lack of punctuality	&	-41.53	^{*}	&	-45.16	^{**}	&	-47.70	^{**}	&	-47.87		^{**}	&	-57.15		^{***}			\\
			Mean ESCS	&	41.93	^{***}	&	36.73	^{***}	&	32.41	^{***}	&	28.23		^{***}	&	25.49		^{**}			\\
			School type (ref. Lyceums)     &  &  & &  &\\																				
			Technical and Vocational 	&	-27.54	^{***}	&	-32.69	^{***}	&	-35.85	^{***}	&	-37.85		^{***}	&	-38.18		^{***}			\\
			Others	&	-90.07	^{***}	&	-106.40	^{***}	&	-113.44	^{***}	&	-106.16		^{***}	&	-95.73		^{***}			\\
			South Italy	&	-39.62	^{***}	&	-43.06	^{***}	&	-43.21	^{***}	&	-41.96		^{***}	&	-42.52		^{***}			\\
			\bottomrule 
			
			\multicolumn{6}{l}{\dag\footnotesize{Point estimates with associated \textit{p}-value: $^{***}\ p<0.01$; $^{**}\ p<0.05$; $^{*}\ p<0.1$}}\\
		\end{tabular}

\end{table}

Findings particularly show that gender has a negative effect on mathematics score at all levels with disadvantage increasing as the mathematics score increases, hence also controlling for other individual and school-related characteristics males still outperform females in math, and the gender gap in math is higher at the upper tail compared with the lower tail of the distribution (the difference between male and female students at the first decile is 14.38, while at the quantile 0.9 is 16.95). The peak of the gender gap is at $q = 0.75$, where boys score about 18 points above girls. The estimated `compositional effect' of gender (i.e. how much the percentage of females in a school affects the mathematics scores in that school), is significantly different from zero along the entire conditional distribution and in addition it increases as the score increases. Namely, if a student moves to a school where the share of female is 10\% higher, his/her score decreases about 5 points for $q = 0.1$ and about 7 points for $q = 0.9$. 

Looking at the other covariates further interesting considerations can be made. As expected, native students perform better than migrants. 
The raw penalty for immigrant students is lower at the bottom of the distribution than elsewhere (the estimated coefficients range from $-14.58$ to $-18.74$), indicating a more considerable negative effect for migrant students at the top end of the distribution.
It may be the result of the difficult for immigrants to achieve very high scores. The share of immigrants in the school shows a negative effect. The effect considerably decreases along the conditional distribution of the mathematics achievement highlighting that high-performing students are less influenced by the presence of immigrants in the school. 
As we expected, math scores of students that never failed school are higher than those who did and the effect of such covariate decreases with increasing quantiles. By contrast, being late for school regularly has a higher effect in the upper tail of distribution, thus the gap is more elevated between best-performing students. If we look at the compositional effects of these two last covariates a further consideration can be made.  
The negative impact of the mean of students who repeat a grade is evident across the entire distribution; moreover this effect increases, in absolute term, as scores increase (the estimated effect is equal to $-51.42$ and $-67.51$ for $q = 0.1$ and $q = 0.9$, respectively). Also the effect of the percentage of students being late at school follows a similar pattern ($-41.53$ and $-57.15$ for $q = 0.1$ and  $q = 0.9$, respectively). Insofar as these two variables are good proxies for the share of truants, findings highlight how the school environment can affect negatively the best performing students. Namely, the presence of negligent students can lead teachers to adopt teaching methods targeted at these more troublesome students to the detriment of good students being less stimulated and consequentially likely affected negatively in their performance.

Furthermore, also the index of anxiety has a greater impact on the upper quantiles, confirming that high-achieving students may be more worried than low-achievers about getting poor grades (see among others \citealp{oecd2016volIII,foley2017}).

Looking at the individual effect of  ESCS index, we found that it is significantly associated with math scores only at $q = 0.90$ ($\beta_{0.9} = 6.07$), whereas its school compositional effect (e.g. the mean of ESCS) is always significantly different from zero. Thus, the effect of peer ``socio-economic status'' tends to have a stronger impact on student performance than the corresponding individual socio-economic status.

As a consequence, being enrolled in an `elite school', i.e. a school with many students from wealthy backgrounds, improves the math scores. Moreover, the effect of the socio-economic context of the school tends to decrease as quantiles increase:  low performing students therefore can enhance their performance in `elite schools'. This last consideration leads us to claim that it is necessary to implement targeted policy measures in order to avoid the ``ghettoisation'' of schools, that is schools with a high concentration of poor students, so that students can benefit from the school-environment and improve their scores.

Finally, both macro-region and type of school have the expected negative sign and their effects are significantly different from zero in the overall distribution. Students living in Southern regions have low math scores than those living in Northern regions, geographical disparities are less evident at the lower tail of the score distribution and then increase even if they are almost stationary as quantiles increase. Students attending technical, vocational and other schools perform worse than those attending Lyceums confirming findings of other studies (see among others \citealp{bratti2007, matteucci2014}).

\section{Simulation study} \label{simulation}
A small Monte Carlo simulation study was carried out to evaluate the performance of the Weighted-MQRE 
regression at three quantiles, $q=0.10, 0.25, 0.50$. The focus of this simulation study is to compare the estimation of the fixed effects of the weighted version of the MQRE and the unweighted version. Moreover, we assess the approximations of the standard errors of the fixed effects. For both aims, finite population values $y_{ij}$ are generated under the two-level random intercept model:
\begin{equation*}
y_{ij}=100+2x_{ij}+\gamma_j + \varepsilon_{ij},  \quad \qquad i=1,...,N_j, \quad j=1,...M.
\end{equation*}

The number of level 2 units in the population is $M=170$, each with the same number of level 1 units, $N_j=50$. The auxiliary variable is uniformly distributed in [0, 20]. 

The level 1 and level 2 residuals are independently generated as follows: $\gamma \sim 0.9\,N(0,1)+0.1\,N(9,20)$ and $\varepsilon \sim 0.9\,N(0,3.3)+0.1\,N(10,75)$. This represents a situation under outlier contamination in both hierarchical levels.

Once the finite population values were obtained, we adopted the following two stage sampling design. Level 2 units are divided into three strata according to whether $\gamma<-1$, $-1\leq \gamma \leq 1$ or $\gamma>1$ and simple random samples of size $0.15\,m$, $0.65\,m$ and $0.20\,m$ are selected from respective strata. The overall number of sampled cluster is $m=100$. 
From each level 2 units, level 1 units are partitioned into two clusters according to whether $\varepsilon > 0 $ or $\varepsilon \leq 0 $ and, also in this case, simple random samples of size $0.75\,n_j$ and $0.25\,n_j$ are selected form respective strata. The size $n_j$ was choice to be proportional to $N_j$, $n_j=0.3\,Nj$. This leads to a total sample size of $n=1500$. This sampling scheme is very similar to the method that was used in \cite{rhs2006} except that we sampled both level units in a different way and in particular we have sampled more cluster and units with extreme values. With this choice, we assure that outliers in the population are also represented in the sample and this allows us to evaluate the performance of the models in a situation where a robust approach is required. 

According to the previous studies, by making the sampling probabilities at level 1 and 2 dependent on the corresponding residuals, we ensure that the sampling scheme is informative at both levels. 

We replicate this scenario $R=500$ times and we compare the fixed effects of the MQRE-2L \citep{tzavidis2016} and the Weighted-MQRE, with tuning constant of the Huber influence function $c=1.345$. At $q=0.5$ we have estimated also the Linear Mixed Model (Eq. \ref{LMM_ri}). Given the informativeness of the sampling probabilities we expect that the Weighted-MQRE performs better than the MQRE. 

For each regression parameter, performance is evaluated using the so called ``Average Relative Bias'' (ARB), defined as:

\begin{equation*}
\mathrm{ARB}(\hat{\theta}) = R^{-1} \sum_{r=1}^{R}\frac{\hat{\theta}^{(r)} - \theta}{\theta} \times 100,
\end{equation*}
where $\hat{\theta}^{r}$ is the estimated parameter at quantile $q$ for the $r^{th}$ replication and $\theta$ is the corresponding `true' value of this parameter.

Table \ref{simB} reports the simulation results for estimators of the fixed effects for $q =$ 0.10, 0.25, 0.50. For the Weighted-MQRE we report results with ``scaling method 2''; estimates using unscaled weights and ``scaling method 1'' are not reported because they are almost identical. Scaling tends to affect the estimation of the variance components, but evaluating the estimation of the variance parameters is not a focus concerning this work (for more detail about the effect of scaling see simulation results in \citealp{pfeffermann1998weights}).
Taking a closer look at ARB for the fixed effects, for the slope we observe that in all models there is almost no bias. However the use of weights have an impact on the estimation of the intercept. Indeed, the bias in the Weighted-MQRE is lower than that of the unweighted models. In particular, the bias tends to decrease with the quantiles. The ARB for intercept of the Weighted-MQRE is around 0.8\% for $q=0.1$ and around 0.6\% for $q=0.5$, while for the MQRE it is around 1.4 for the first decile and around 1.1\% for $q=0.5$. The LMM appears to be the worst model, the bias at $q=0.5$ is around 2\%; given the presence of the outliers in both levels, also the weighted version of the LMM shows a severe bias (the results are not reported here). 

\begin{landscape} 
\begin{table}\caption{\label{simB}Values of bias ARB and the average of point estimates over simulations of the fixed effects using Weighted-MQRE, MQRE and LMM at $q=0.10,0.25,0.50$} 
	\centering
	\small
	
	\begin{tabular} {lrrrrrrrrrrrr}
		\toprule 
		\noalign{\medskip}
		\textit{Method}	& \multicolumn{2}{c}{$\hat{\beta_0}$,$q=0.10$}	& \multicolumn{2}{c}{$\hat{\beta_1}$,$q=0.10$} & \multicolumn{2}{c}{$\hat{\beta_0}$,$q=0.25$} & \multicolumn{2}{c}{$\hat{\beta_1}$,$q=0.25$} & \multicolumn{2}{c}{$\hat{\beta_0}$,$q=0.50$} & \multicolumn{2}{c}{$\hat{\beta_1}$,$q=0.50$}\\
		\noalign{\medskip}
		\cmidrule(lr){2-3} \cmidrule(lr){4-5} \cmidrule(lr){6-7} \cmidrule(lr){8-9} \cmidrule(lr){10-11} \cmidrule(lr){12-13}
		\noalign{\smallskip}
		& \multicolumn{1}{c}{\textit{ARB}} & \multicolumn{1}{c} {$\hat{\beta_0}$ } & \multicolumn{1}{c}{\textit{ARB}}   &  \multicolumn{1}{c} {$\hat{\beta_1}$ }   & \multicolumn{1}{c}{\textit{ARB}}   &   \multicolumn{1}{c} {$\hat{\beta_0}$ } & \multicolumn{1}{c}{\textit{ARB}}   &   \multicolumn{1}{c} {$\hat{\beta_1}$ } & ARB  &   \multicolumn{1}{c} {$\hat{\beta_0}$ } & \multicolumn{1}{c}{\textit{ARB}}   &  \multicolumn{1}{c} {$\hat{\beta_1}$ }  \\
		\noalign{\smallskip}
		\midrule\noalign{\medskip}
		Weighted-MQRE & 0.797  & 98.452 & -0.025  & 1.9995  & 0.675  & 99.443 & -0.013  & 1.9997& 0.626  & 100.627 & -0.011 & 1.9998 \\
		
		\noalign{\smallskip}
		
		\noalign{\medskip}
		MQRE & 1.406  & 99.048 & -0.025  & 1.9995  & 1.258 & 100.020 & -0.016 & 1.9997 & 1.118 & 101.121 & -0.013 & 1.9997 \\
		
		\noalign{\smallskip}
		\noalign{\medskip}
		LMM & - & - & - & - & - & -  & - & -   & 2.022  & 102.026 & -0.013  & 1.9997  \\
		\noalign{\smallskip}
		\bottomrule	
		
	\end{tabular}
\end{table}
\end{landscape}

\begin{table}\caption{\label{seSim}Empirical and estimated standard Errors of he fixed effects for $q=0.10,0.25,0.50$ using Weighted-MQRE with tuning constant $c=1.345$} 
	\centering
	\small
		\begin{tabular} {lcccc}
		\toprule 
		\textit{	Values of q} & \multicolumn{2}{c}{$\hat{\beta_0}$}	& \multicolumn{2}{c}{$\hat{\beta_1}$} \\
		\noalign{\smallskip}
		\cmidrule(lr){2-3} \cmidrule(lr){4-5}
		\noalign{\smallskip}
		& \multirow{1}{*}{\textit{Empirical}} & \multirow{1}{*}{\textit{Estimated}} & \multirow{1}{*}{\textit{Empirical}}
		& \multirow{1}{*}{\textit{Estimated}}       	 \\
		& \multirow{1}{*}{\textit{standard}} & \multirow{1}{*}{\textit{standard}} & \multirow{1}{*}{\textit{standard}}
		& \multirow{1}{*}{\textit{standard}}       	 \\
		& \multirow{1}{*}{\textit{error}} & \multirow{1}{*}{\textit{error}} & \multirow{1}{*}{\textit{error}}
		& \multirow{1}{*}{\textit{error}}       	 \\
		\noalign{\smallskip}
		\midrule 
		\noalign{\smallskip}
		$q=0.10$ & 0.198 & 0.211 & 0.016 & 0.015 \\
		$q=0.25$ & 0.169 & 0.190 & 0.013 & 0.012 \\
		$q=0.50$ & 0.180 & 0.215 &0.012  & 0.012 \\
		\bottomrule				
	\end{tabular}
\end{table}
Having evaluated the performance of the Weighted-MQRE regression, we want to test the performance of the sandwich estimator for the variability of the fixed effects. Therefore, we compare the empirical standard errors and estimated standard errors. For each estimator $\hat{\theta}$, at $q=$0.10, 0.25, 0.50, Table \ref{seSim} reports averages over simulations of the Monte Carlo standard error $S(\hat{\theta})=\sqrt{ \bigl[ R^{-1} \sum_{r=1}^{R}( \hat{\theta}^{(r)}- \bar{\theta})^2 \bigr]}$, where $\bar{\theta}= R^{-1} \sum_{r=1}^{R}\hat{\theta}^{(r)}$, and the estimated standard errors of the fixed effects $\hat{\beta}_q$.
It can be observed that the estimated standard error of the estimators at all values of $q$ offers a good approximation to the empirical standard errors. Furthermore, as we expected, inserting weights into the estimation procedure it produces a little overestimation of the standard errors of the intercept.  

In summary, our simulation results showed that accommodating sampling weights in the estimation procedure of the MQRE is a good strategy for dealing with an informative sampling design and outlier contamination in the data.

\section{Conclusion} \label{conclusion}

In this paper, we focus on the gender effect on math outcomes at upper secondary school in Italy. Previous research has suggested that the female penalty may help explain the  low share of women in STEM education, which is in turn one of the reasons of the poor female presence in research and innovation decision-making positions.
We improve existing studies in several ways but, most importantly, we propose a new methodology that allows to apply M-quantile regression in a multilevel framework and in complex survey design. We also propose a simulation study in order to evaluate the performance of the fixed effects of the Weighted-MQRE showing that the weighted estimator can reduce the bias induced by the complex sampling design when subjects and clusters are selected with probabilities that depend on the model's random terms.
To do this we use data drawn from 2015 PISA survey for Italy. 

In summary, the existence of a gender-related issue in the Italian education system persists not only around the mean but in the overall distribution and, in particular, the ability of girls in school is not able to mitigate the effect of being female. Indeed, as we pointed out in the Section 1, high-performing students tend to continue studying and, moreover, students who obtain high grades in mathematics tend to choose science-related subjects in tertiary education \citep{oecd12gender}. As a consequence, we can argue that the larger gender gap among the top performers is one of the main reasons of the low share of girls in STEM. Reducing this gap is essential to combat gender-segregation in the tertiary education and this might in turn reduce the gender gap in labour market outcomes. 
The above findings can be linked with some interesting findings of the \cite{abc2015} report. Student attitudes and beliefs are correlated with educational achievement, and in particular, the lack of self-confidence in mathematics among girls affects their mathematics' achievement. Indeed, self-confidence gives students the possibility to go through a trial-and-error process that is very important in developing mathematical skills; moreover, even the high-achieving girls have low levels of confidence in their ability and this can explain the gender gap in mathematics. Accordingly, encouraging girls to change their attitudes could improve gender equality in the math scores, and parents and teachers represent the best mediators in this process, because they directly affect their education and attitudes. 

It is beyond doubt that psychological factors frame much of our tools for reducing these gaps, but the macro-societal characteristics (i.e. societal inequality and the welfare orientation) of Italy may play an important role as well \citep{marks2008, cipollone2014}.  
In this view, the welfare state of Italy, combining elements of the continental-corporatist model and the familist model model, is characterized by relatively low female labour market participation anyhow \citep{ferrera2010, katrougalos2008, naldini2013, lynch2014}. 
Therefore, policy innovations (family-oriented policies) are likely another important way of reducing gender gaps.

Finally, further work using other OECD countries, and comparing results across countries would be the natural extension for applying the new methodology and analysing cross-country level effects that can add further information for improving gender equality in Italy.  

\section*{Appendix A} 
The steps of the estimation algorithm are as follows:
\begin{itemize} 
	\item[1.] Let $\suq$ and $\seq$ be known.
	\item[2.] Given the variance parameter and the covariance matrix $\VV_{q}$, the iterative equation for the estimates of $\BB_{\psi q}$ is:
	\begin{equation*} 
	\BB_{\psi q}^{t+1}=\BB_{\psi q}^{t}+ [\XX^{T} \Gam \UU_q^{-1/2} \DD_q(\BB_{\psi q}^{t}) \UU_q^{1/2} \VV_q^{-1} \XX]^{-1} [\XX^{T} \Gam \VV_q^{-1} \UU_q^{1/2} \psi_q(\rr_q)],
	\end{equation*}
	
	where $\Gam$ is a diagonal matrix with $j${th} diagonal element equals to the second level weight $w_j$; $\DD_q(\BB_{\psi q}^{t})$ is a diagonal matrix with  $j${th} diagonal element  $D_{ijq}=\psi_q^{'}(r_{ijq})=\partial\psi_q^{'}(r_{ijq})/\partial r_{ijq}$.
	
	\item[3.] The estimates of $\BB_{\psi q}$ are used to obtain the estimates of the variance parameters with a fixed point iterative method. As pointed out in \cite{tzavidis2016} this requires the change of the estimating equation in Eq. \ref{var} to:
	
	\begin{equation*} 
	\begin{split}
	& \Gam \psi_q(\rr_q)^{T} \UU^{1/2}_q \VV^{-1}_q \ZZ\ZZ^{T} \VV^{-1}_q \UU^{1/2}_q \psi_q(\rr_q) \\ &
	- K_{2q}\mathrm{tr}\left\{\Gam \VV^{-1}_q \ZZ\ZZ^{T}\ \VV^{-1}_q (\ZZ\ZZ^{T}\, \WW^{-1}) \begin{pmatrix}
	\suq\\
	\seq 
	\end{pmatrix} \right\} = 0
	\end{split}
	\end{equation*}
	
	\begin{equation*}
	\begin{split} 
	& \Gam \psi_q(\rr_q)^{T} \UU^{1/2}_q \VV^{-1}_q \WW^{-1} \VV_q^{-1} \UU^{1/2}_q \psi_q(\rr_q) \\ & 
	- K_{2q}\mathrm{tr}\left\{\Gam \VV^{-1}_q \WW^{-1}\ \VV^{-1}_q (\ZZ\ZZ^{T}\, \WW^{-1}) \begin{pmatrix}
	\suq \\
	\seq
	\end{pmatrix} \right\} = 0
	\end{split}
	\end{equation*}

	The fixed point algorithm of the estimating equation for the $t^{th}$ iteration is the following:
	\[
	\begin{pmatrix}
	\suq \\
	\seq
	\end{pmatrix}=
	\left\{ \mathbf{A} \begin{pmatrix}  \suqt \\
	\seqt \end{pmatrix} \right\}^{-1}   a \begin{pmatrix}  \suqt \\
	\seqt \end{pmatrix},
	\]
	
	where 
	
	\[
	\mathbf{A}\begin{pmatrix}
	\suq \\
	\seq
	\end{pmatrix}=\begin{pmatrix} 
	K_{2q} \mathrm{tr}(\Gam \VV_q^{-1} \ZZ\ZZ^{T}\VV_q^{-1} \ZZ\ZZ^{t}) & K_{2q} \mathrm{tr}(\Gam \VV_q^{-1} \ZZ\ZZ^{T}\VV_q^{-1} \WW^{-1}) \\
	\\
	K_{2q} \mathrm{tr}(\Gam \VV_q^{-1} \WW^{-1} \VV_q^{-1} \ZZ\ZZ^{T}) & K_{2q} \mathrm{tr}(\Gam \VV_q^{-1} \WW^{-1} \VV_q^{-1} \WW^{-1})
	\end{pmatrix}
	\]

	and
	
	\[
	a\begin{pmatrix}
	\suq \\
	\seq 
	\end{pmatrix}=\begin{pmatrix} 
	\frac{1}{2} \psi_q(\rr_q)^{T} \UU_q^{1/2} \VV_q^{-1} \ZZ\ZZ^{T} \Gam \VV_q^{-1} \UU_q^{1/2} \psi_q(\rr_q)\\
	\\
	\frac{1}{2} \psi_q(\rr_q)^{T} \UU_q^{1/2} \VV_q^{-1} \WW^{-1} \Gam \VV_q^{-1} \UU_q^{1/2} \psi_q(\rr_q)
	\end{pmatrix}
	\]
	
	\item[4.] Iterate steps 2 and 3 until convergence.
\end{itemize}

\section*{Appendix B}
The asymptotic covariance matrix of the estimators can be written as $\mathrm{G}^{-1}\mathrm{F}(\mathrm{G}^{-1})^T$. 

The weighted estimating equations in the Eq. \ref{beta} and Eq. \ref{var} can be written as:

\begin{equation} \label{B1}
\sum_{j=1}^{m} w_j \boldsymbol \Phi_{qj}(\boldsymbol \theta_q)=\mathbf{0},
\end{equation}

where $\boldsymbol \theta_q=(\hat{\BB}_q^T,\hat{\sigma}^2_{\varepsilon_q}, \hat{\sigma}^2_{\gamma_q})$ and $\boldsymbol \Phi_{qj}(\boldsymbol \theta_q)= \big(\boldsymbol \Phi_{qj\BB_q}^T,  \Phi_{qj \suq},  \Phi_{qj \seq} \big)^T$. 

Under a general response distribution $D$, the estimator $\hat{\boldsymbol\theta}_q$ satisfying equations (\ref{B1}) is estimating a root $\boldsymbol \theta_q$ of:
$$
\sum_{j=1}^{m}E_D[ w_j \boldsymbol \Phi_{qj}(\boldsymbol \theta_q)]=\mathbf{0}.
$$

where $\boldsymbol \Phi_{qj}(\boldsymbol \theta_q)= \big(\boldsymbol \Phi_{qj\BB_q}^T,  \Phi_{qj \suq},  \Phi_{qj \seq} \big)^T$ for particular choice of $\boldsymbol \Phi_{qj}(\boldsymbol \theta_q)$. Under a general response distribution of $D$ (i.e. not necessarily Gaussian), a robust estimator is estimating a root $\boldsymbol \theta_q$ of the expectation with respect to $D$ of the previous estimating equations. 

Provided that 
$$ n^{-1} \sum_{j=1}^{m} E_D \Bigg[ \dfrac{\partial (w_j \boldsymbol \Phi_{qj}(\boldsymbol \theta_q))}{\partial \boldsymbol \theta_q} \Bigg]  \rightarrow \mathbf{G}, $$

is the diagonal block information matrix of dimension $(p+2) \times (p+2)$ and it is positive definite, and 
$$ n^{-1} \sum_{j=1}^{m} E_D \Bigg[ \dfrac{(w_j \boldsymbol \Phi_{qj}(\boldsymbol \theta_q))^T}{w_j \boldsymbol \Phi_{qj}(\boldsymbol \theta_q)} \Bigg]\rightarrow \mathbf{F}, $$
is the the matrix of the variance of the normalised score function, a Taylor series approximation which which holds uniformly in a neighbourhood of $\boldsymbol \theta_q$ 
is: $$ \hat{\boldsymbol \theta}_q \approx \boldsymbol \theta_q + \mathbf{G}^{-1}n^{-1} \sum_{j=1}^{d} w_j \boldsymbol \Phi_{qj}(\boldsymbol \theta_q) + o_p(n^{-1/2}), \quad n\rightarrow \infty.$$

Following \cite{richardsonwelsh1995} the covariance matrix can be consistently estimated by $\mathbf{\hat{G}}^{-1}\mathbf{\hat{F}}({\mathbf{\hat{G}}^{-1}})^T$ where  the matrices $\mathbf{\hat{G}}$ and $\mathbf{\hat{F}}$ are evaluated at $\hat{\boldsymbol \theta_q}$.

\section*{Appendix C}
Following \cite{fabrizi2014w} we rewrite the assumption 1, 2, 8 and 9 of \cite{wang2011}. 

\textbf{Assumption A1.} The expected number of cluster
$m^{*}=E_d(m|\mathbf{y}_1,\mathbf{y}_2,...,\mathbf{y}_M)=O(M^\delta)$, where the vectors $\mathbf{y}_i$, $i=1,2,...3,M$ denote all the variables of interest in each group of the $M$th population. The subscript $d$ denotes the expectations with respect to the randomization distribution induced by the sampling design.

\textbf{Assumption A2.} $K_L\leq M \pi_j/m^* \leq K_U$ for all j, with $K_L$ and $K_U$ positive constants.

\textbf{Assumption A3.} For any vector $\mathbf{z}$ with finite $2+\lambda$ moments with arbitrarily small $\lambda > 0$, we assume $V_d(\bar{\mathbf{z}}_\pi|\mathbf{y}_2,...,\mathbf{y}_M) \leq c_1 m^{*-1}(M-1)^{-1} \sum_{j=1}^{M} (\mathbf{z}_j -\mathbf{\bar{z}}_M)(\mathbf{z}_j -\mathbf{\bar{z}}_M)^T$, for some constant $c_1$ and $\bar{\mathbf{z}}_\pi=M^{-1} \sum_{j=1}^{m} w_j \mathbf{z}_j$.

\textbf{Assumption A4.} The population parameter $\mathbf{B}_{\psi q}$ lies in a closed interval $\Theta_{\mathbf{B}_{\psi q}}$ on $\mathbb{R}^p$.

\textbf{Assumption A5.} The population estimating function $S_M(\cdot)$ and the function $\psi(\cdot)$ satisfy:
\begin{itemize}
	\item[-] the function $\psi(\cdot)$ is bounded;
	\item[-] the population estimating function $S_M(\BB_{\psi q})$ converges to $S(\BB_{\psi q})$ uniformly on $\Theta_{\mathbf{B}_{\psi q}}$ as $M\rightarrow \infty$ and the equation $S(\BB_{\psi q})=0$ has a unique root in the interior of  $\Theta_{\mathbf{B}_{\psi q}}$;
	\item[-] the limiting function $S(\BB_{\psi q})$ is strictly increasing and absolutely continuous with the finite first derivative in $\Theta_{\mathbf{B}_{\psi q}}$, and the derivative $S^{'}(\BB_{\psi q})$ is bounded away from 0 for $\BB_{\psi q}$ in $\Theta_{\mathbf{B}_{\psi q}}$;
	\item[-] the population quantities:
	$$
	\underset{\BB_{\psi q} \in \Theta_{\BB_{\psi q}}}{\mathrm{sup}} M^{\alpha}|S_M(\mathbf{B}_{\psi q} + M^{-\alpha}\BB_{\psi q}) - S_M(\mathbf{B}_{\psi q}) - S(\mathbf{B}_{\psi q} + M^{-\alpha}\BB_{\psi q}) + S(\mathbf{B}_{\psi q})| \rightarrow 0
	$$
	and 	
	\begin{equation*}
	\begin{split}
	\underset{\BB_{\psi q} \in \Theta_{\BB_{\psi q}}}{\mathrm{sup}} & M^{-1} \sum_{j=1}^{M} |[ \psi_q(\VV_j^{-1} (\mathbf{y}_j  - \xx_j^T \BB_{\psi q} - M^{-\alpha}\xx_j^T \BB_{\psi q}))  - \\ & \psi_q(\VV_j^{-1} (\mathbf{y}_j  - \xx_j^T \mathbf{B}_{\psi q}))] w_j \XX_j \VV_j^{-1/2}| = O_p(M^{-\alpha})
	\end{split}
	\end{equation*}
	where $\Theta_{\BB_{\psi q}}$ is a large enough compact set in $\mathbb{R}^p$ and $\alpha \in (1/4, 1/2)$.
\end{itemize}

\bibliographystyle{chicago} 
\bibliography{PISAbib,MQRE}	

\end{document}